\documentclass{conm-p-l}

\usepackage{amsgen,amsmath,amstext,amsbsy,amsopn,amsfonts,amssymb, enumerate}

\newcounter{mycounter}

 \newcounter{mycounter2}

\def\~{\tilde}

\def\F{{\mathbb F}}

\def\Zi{{\mathbb Z}}

\def\cFp{{\cal F}_p}

\def\w{{\omega}}

\def\dim{{\rm dim}}

\def\gr{{\rm gr}}

\def\cFp{{\mathcal F}_p}

\usepackage{amssymb}

\newtheorem{defi}{Definition}[section]
\newtheorem{thm}[defi]{Theorem}
\newtheorem{lem}[defi]{Lemma}
\newtheorem{cor}[defi]{Corollary}

\newtheorem{eg}[defi]{Example}
\newtheorem{prop}[defi]{Proposition}
\newtheorem{rem}[defi]{Remark}

\DeclareMathOperator{\la}{\langle}
\DeclareMathOperator{\ra}{\rangle}
\DeclareMathOperator{\su}{\subseteq}
\DeclareMathOperator{\charac}{char}
\DeclareMathOperator{\ad}{ad}

\newcommand{\uea}[1]{U(#1)}

\begin{document}

\title[Isomorphism invariants of  enveloping algebras]{Isomorphism invariants of  enveloping algebras}

\author{\textsc{Hamid Usefi}}
\address{Department of Mathematics and Statistics,
Memorial University of Newfoundland,
St. John's, NL,
Canada, 
A1C 5S7}
\email{usefi@mun.ca}

\thanks{The research of the author was supported by NSERC of Canada.}

\begin{abstract}
Let $L$ be a Lie algebra with its enveloping algebra $U(L)$ over a field.
In this paper we survey results concerning the isomorphism problem for enveloping algebras: given another Lie algebra $H$ for which
$U(L)$ and $U(H)$ are isomorphic as associative algebras, can we deduce that $L$ and $H$ are isomorphic Lie algebras?
Over a field of positive characteristic we consider a similar problem for restricted Lie algebras, that is, given restricted Lie algebras 
$L$ and $H$ for which their restricted enveloping algebras  are isomorphic as algebras, can we deduce that
$L$ and $H$ are isomorphic?
\end{abstract}
\date{\today}


\date{\today}

\maketitle

\section{Introduction}
Let $L$ be a Lie algebra with universal enveloping algebra $U(L)$ over a filed $\F$.
Our aim in this paper is to survey the results concerning the
isomorphism problem for enveloping algebras: given another Lie algebra 
$H$ for which $U(L)$ and $U(H)$ are isomorphic as associative algebras, can we deduce that 
$L$ and $H$ are isomorphic Lie algebras? We can ask weaker questions  in the sense that
given $U(L)\cong U(H)$, what invariants of $L$ and $H$ are the same?
 We  say that a particular invariant of $L$ is
{\em determined} (by $U(L)$), if
every Lie algebra $H$ also possesses this invariant whenever
$U(L)$ and $U(H)$ are isomorphic as associative algebras. For example, it is well-known that the dimension of a
finite-dimensional Lie algebra $L$ is determined by $U(L)$ since
it coincides with the Gelfand-Kirillov dimension of $U(L)$. 

 The closely related isomorphism problem for group rings asks: is every finite group $G$
determined by its integral
group ring ${\mathbb Z}G$?  A positive solution for the class of all
nilpotent groups was given independently in \cite{RS} and
\cite{W}. There exist, however, a pair of non-isomorphic finite solvable groups of derived length 4 whose integral group rings are isomorphic (see \cite{H}). 

In Section \ref{fox}, we discuss  identifications of certain Lie subalgebras associated to the augmentation ideal $\w(L)$ of $U(L)$. Let  $S$ be a Lie subalgebra of $L$.
 The identification of the subalgebras $L\cap \w^n(L)\w^m(S)$ naturally
arises in the context of enveloping algebras.

 In Section \ref{aug}, we show that certain invariats of $L$ are determined by $U(L)$ including the nilpotence class of a nilpotent Lie algebra $L$. The results of this section motivate one to investigate the isomorphism problem in detail for low dimensional nilpotent Lie algebras. The isomorphism problem for nilpotent Lie algebras of dimension at most 6 is discussed in Section \ref{low-dim}. It turns out that there exist counterexamples to the isomorphism problem  in dimension 5 over a field of characteristic 2 and in dimension 6 over a field of characteristic 2 and 3.
The conclusion is that as the dimension of $L$ increases we have to exclude more  fields of  positive characteristic to avoid counterexamples. Indeed, if there is a pair of Lie algebras $L$ and $H$ that provides a counterexample over a field $\F$ then a central extension of $L$ and $H$ would provide a counterexample in a higher dimension over $\F$.  Furthermore, it is shown    in Example \ref{ex1} that 
over any field of positive characteristic $p$ there exist non-isomorphic Lie algebras of  dimension $p+3$ whose enveloping algebras are isomorphic. So over any field of positive characteristic $p$ and any integer  $n\geq p+3$ there exists a pair of non-isomorphic Lie algebras of dimension $n$ whose enveloping algebras are isomorphic.

These observations convince us to consider the isomorphism problem over a field of characteristic zero, however there are other invariants that we expect to be determined over any filed. Some of these questions are listed in Section \ref{quo}.

Nevertheless, it makes sense to consider the isomorphism problem for restricted Lie algebras over a filed $\F$ of positive characteristic $p$. We denote the restricted enveloping algebra of a restricted Lie algebra $L$ by $u(L)$. Given another restricted Lie algebra $H$ for which $u(L)\cong u(H)$ as associative algebras, we ask what invariants of $L$ and $H$ are the same? For example, since $\dim_{\F} u(L)= p^{\dim_{\F} L}$, the dimension of $L$ is determined. Unlike abelian Lie algebras whose only invariant
is their dimension, abelian restricted Lie algebras have more structure. As a first step,   abelian restricted Lie algebras  are considered in Section \ref{rest}. Other known results for restricted Lie algebras are also discussed in Section \ref{rest}.

In Section \ref{hopf}, we have collected the known results about the Hopf algebra structure of $U(L)$ and $u(L)$ deducing  that 
the isomorphism problem is trivial if the Hopf algebra structures of $U(L)$ or $u(L)$ are considered. In this section an example is
given showing that the enveloping algebra of a Lie superalgebra $L$ may not necessary determine the dimension of $L$.
 As we mentioned earlier some open problems are discussed in Section \ref{quo}.

\section{Preliminaries}

In this section we collect some basic definitions that can be   found in  \cite{B} or 
\cite{SF}. Every associative algebra can be viewed as a Lie algebra under the natural Lie bracket
$[x,y]=xy-xy$. In fact, every Lie algebra can be embedded into an associative algebra
in a canonical way.

\begin{defi} 
Let $L$ be a Lie algebra over $\F$ and $U(L)$ an associative algebra.
Let $\iota:L\to U(L)$ be a Lie homomorphism. The pair $(U(L),\iota)$ is called a (universal) 
enveloping algebra of $L$ if for every associative algebra $A$ and every Lie homomorphism 
$f:L\to A$ there is a unique algebra homomorphism ${\bar f}:U(L)\to A$ such that
${\bar f}\iota=f$.
\end{defi}
 It is clear that if an enveloping algebra exists, then it is unique up to isomorphism. Its 
existence can be shown as follows. Let $T(L)$ be the tensor algebra based
 on the vector space $L$, that is
 $$
 T(L)=\F\oplus L \oplus (L\otimes L)\oplus\cdots.
 $$
The multiplication in $T(L)$ is induced by concatenation which turns $T(L)$ into an
associative algebra. Now let $I$ be the ideal of $T(L)$ generated by all
elements of the form
$$
[x,y]-x\otimes y+y\otimes x, ~~x,y\in L,
$$
and let $U(L)=T(L)/I$. If we denote by $\iota$ the restriction to $L$ of the
natural homomorphism $T(L)\to U(L)$, then it can be verified that  
$(U(L),\iota)$ is the enveloping algebra of $L$. Furthermore, since $\iota$
is injective, we can regard $L$ as a Lie subalgebra of $U(L)$. In fact we can say more:

\begin{thm}[Poincar\' e-Birkhoff-Witt]
Let $\{x_j\}_{j\in {\mathcal J}}$ be a totally-ordered basis for $L$ over $\F$.
Then $U(L)$ has a basis consisting of PBW monomials, that is, monomials  of the form
$$
x_{j_1}^{a_1}\cdots x_{j_t}^{a_t},
$$
where $j_1<\cdots< j_t$ are in ${\mathcal J}$ and $t$ and each $a_i$ are non-negative 
integers.
\end{thm}

This result is commonly referred to as the PBW Theorem.
Let $H$ be a subalgebra of $L$. It follows from the PBW Theorem that the extension of 
the Lie homomorphism $H\hookrightarrow L\hookrightarrow U(L)$ to $U(H)$ is an injective 
algebra homomorphism. So we can view  $U(H)$ as a subalgebra of $U(L)$.

Next consider the \emph{augmentation map} $\varepsilon_L: U(L)\to \F$
which is the unique algebra homomorphism induced by $\varepsilon_L(x)=0$ for
every $x\in L$. The kernel of $\varepsilon _L$ is called the \emph{augmentation ideal} of $L$ 
and will be denoted by $\omega(L)$; thus, $\w(L)=LU(L)=U(L)L$.
 We denote by $\w^n(L)$ the
$n$-th power of $\w(L)$ and $\w^0(L)$ is $U(L)$.

We consider left-normed commutators, that is
$$
[x_1, \ldots, x_n]=[[x_1, x_2], x_3, \ldots, x_n].
$$
The \emph{lower central series} of $L$  is defined inductively by 
$\gamma_1(L)=L$ and  $\gamma_{n}(L)=[\gamma_{n-1}(L),L]$. The second term 
will be also denoted by $L'$; that is, $L'=\gamma_2(L)$. If $L'=0$ then $L$ is called 
abelian. A Lie algebra $L$ is said to be \emph{nilpotent} if $\gamma_n(L)=0$ for
some $n$; the \emph{nilpotence class} of $L$ is the minimal integer $c$
such that $\gamma_{c+1}(L)=0$. Also, $L$ is called \emph{metabelian} if $L'$ is abelian.
Let $L$ be a Lie algebra  over a field $\F$ of positive  characteristic $p$ and denote by 
$\ad  : L \to L$ the adjoint representation of $L$ given by
$(\ad x)(y) = [y,x]$, where $x, y \in L$. 
 A mapping $^{[p]}:L\to L$ that satisfies the following 
properties for every $x,y\in L$ and $\alpha\in \F$:

\begin{enumerate}
\item $(\ad x)^p=\ad(x^{[p]})$;
\item $(\alpha x)^{[p]}=\alpha^p x^{[p]}$; and,
\item $(x+y)^{[p]}=x^{[p]}+y^{[p]}+\sum_{i=1}^{p-1} s_i(x,y)$,
where $is_i(x,y)$ is the coefficient of $\lambda^{i-1}$ in $\ad(\lambda x+y)^{p-1}(x)$.
\end{enumerate}
is called a $[p]$-mapping. The pair $(L, [p])$ is called a \emph{restricted Lie algebra}.

\begin{rem}\label{si-gammap}
By expanding $\ad(\lambda x+y)^{p-1}(x)$ it can be seen  that $s_i(x,y)\in \gamma_p(\la 
x,y\ra)$, for every $i$.
\end{rem}

Every Associative algebra  can be regarded as a restricted Lie algebra  with the natural Lie 
bracket and exponentiation by $p$ as the $[p]$-mapping: $x^{[p]}=x^p$.

A Lie subalgbera $H$ of $L$ is called a \emph{restricted subalgebra} of $L$ if $H$ is closed 
under the $[p]$-map.
 
\begin{defi} 
Let $L$ be a restricted Lie algebra over $\F$ and $u(L)$ an associative algebra.
Let $\iota:L\to u(L)$ be a restricted Lie homomorphism. The pair $(u(L),\iota)$ is called a 
restricted (universal) enveloping algebra of $L$ if for every associative algebra $A$ and 
every restricted Lie homomorphism 
$f:L\to A$ there is a unique algebra homomorphism ${\bar f}:u(L)\to A$ such that
${\bar f}\iota=f$.
\end{defi}
It is clear that if  restricted enveloping algebras exist then they are unique up to
an algebra isomorphism. Let $I$ be the ideal of $U(L)$ generated 
by all elements $x^{[p]}-x^p$, $x\in L$. Put $u(L)=U(L)/I$. Then 
$(u(L),\iota)$ has the desired property, where $\iota$ is the restriction to $L$
of the natural map $U(L)\to u(L)$.
  The analogue of the PBW Theorem for restricted Lie algebras is due to Jacobson:

\begin{thm}[Jacobson]
Let $\{x_j\}_{j\in {\mathcal J}}$ be a totally-ordered basis for $L$ over a filed $\F$ of positive characteristic $p$.
Then $u(L)$ has a basis consisting of restricted PBW monomials.
\end{thm} 
An important consequence is that we may regard $L$ as a restricted subalgebra of $u(L)$.
Thus, the $p$-map in $L$ is usually denoted by $x^p$.
 
Let $X$ be a subset of $L$. The restricted subalgebra generated by $X$ in $L$, denoted by 
$\langle X \rangle_p$, is  the smallest restricted subalgebra containing $X$. Also, $X^{p^j}$
denotes the restricted subalgebra generated by all $x^{p^j}$ with $x\in X$. 
Recall that $L$ is $p$-nilpotent if there exists a positive integer $k$ such that $x^{p^k}=0$, for all $x\in L$. We say  $L\in \cFp$ if $L$ is finite dimensional and $p$-nilpotent. Note that if $L\in \cFp$ then $L$ is nilpotent by Engel's Theorem. We denote by $L'_p$ the restricted subalgebra of $L$ generated by $L'$.

\section{Fox-type problems}\label{fox}

Let $L$ be a Lie algebra and $S$ a subalgebra of $L$ over a field $\F$. 
 The identification of the subalgebras $L\cap \w^n(L)\w^m(S)$ naturally
arises in the context of enveloping algebras.  It is proved in \cite{Ri, RU} that $L\cap\w^n(L)=\gamma_n(L)$, for every integer $n\geq 1$.
Furthermore, we have:

\begin{prop}[\cite{RU}]
Let $S$ be a subalgebra of a Lie algebra $L$. The following statements hold for every integer $n\geq 1$.
\begin{enumerate}
\item $\w(S)\cap \w^n(S)\w(L)=\w^{n+1}(S)$; hence, $L\cap \w^n(S)\w(L)=\gamma_{n+1}(S)$ .
\item $\w(S)\cap \w^n(S)U(L)=\w^{n}(S)$.
%\item  $\w^n(S)/\w^{n+1}(S)$ embeds into $\w^n(S)U(L)/\w^n(S)\w(L)$.
%\item  $\w^n(S)/\w^{n+1}(S)$ embeds into $\w^n(S)U(L)/\w^{n+1}(S)U(L)$.
\end{enumerate}
\end{prop}

Hurley and Sehgal in \cite{HSe} proved that if $F$ is  a free group and  $R$  a normal subgroup of $F$, then
$$
F\cap (1+\w^2(F)\w^n(R))=\gamma_{n+2}(R)\gamma_{n+1}(R\cap \gamma_2(F)),
$$
 for every positive integer $n$.  The analogous result for Lie algebras is as follows:

\begin{thm}[\cite{U-JA08}]
Let $L$ be a Lie algebra and $S$ a Lie subalgebra  $L$. For every positive integer $n$, the following subalgebras
of $L$ coincide.
\begin{enumerate}
\item  $\gamma_{n+2}(S)+\gamma_{n+1}(S\cap \gamma_2(L))$,
\item  $L\cap (\w^{n+2}(S)+\w(S\cap \gamma_2(L))\w^n(S))$,
\item  $L\cap \w^2(L)\w^n(S)$.
\end{enumerate}
\end{thm}

The motivation for this sort of problems
also comes from its group ring counterpart.
Let $F$ be a free group, $R$  a normal subgroup of $F$, and
denote by $\mathfrak{r}$ the kernel of the natural homomorphism $\Zi F\to \Zi(F/R)$.
Recall that the augmentation ideal $\mathfrak{f}$ of
the integral group ring $\Zi F$ is the kernel of the map  $\Zi F\to \Zi$ induced by $g\mapsto 1$ for every $g\in
F$.  Fox introduced in \cite{Fox} the problem of identifying the subgroup $F\cap (1+\mathfrak{f}^n\mathfrak{r})$ in terms of $R$.  Following Gupta's initial work on Fox's problem (\cite{Gup}),   Hurley (\cite{Hur2})  and Yunus
(\cite{Yun}) independently gave a complete solution to this problem. 

At the same time Yunus considered the Fox problem for free Lie algebras. Let $\mathcal{L}$ be a free Lie
 algebra and $\mathcal{R}$
an ideal of $\mathcal{L}$. Yunus in \cite{Yun2} identified the subalgebra $\mathcal{L} \cap \w(\mathcal{R})\w^n(\mathcal{L})$ in terms of $\mathcal{R}$. The solutions to the Fox problem for free restricted Lie algebras  is as follows.

\begin{thm} [\cite{U-JA08}]
Let $\mathcal{R}$ be a restricted ideal of a free restricted Lie algebra $\mathcal{L}$. Then 
$$
\mathcal{L}\cap \w^n(\mathcal{L})\w(\mathcal{R})=\sum [\mathcal{R}\cap\gamma_{i_1}(\mathcal{L}),\ldots,\mathcal{R}\cap\gamma_{i_k}(\mathcal{L})]^{p^j}+\sum 
(\mathcal{R}\cap\gamma_i(\mathcal{L}))^{p^{\ell}},
$$
 where the first sum is over all tuples  $(i_1,\ldots,i_k)$, $k\geq 2$,  and non-negative integer $j$ such that 
$p^j(i_1+\cdots+i_k)-i_t\geq n$, for every $t$ in the range $1\leq t\leq k$ and the second sum is over all  
positive integers $i$ and $\ell$ such that $(p^{\ell}-1)i\geq n$.
 \end{thm}

 Let $L$ be a restricted Lie algebra.
The dimension subalgebras of $L$ are defined by
$$
D_n(L)=L\cap \w^n(L).
$$
\begin{thm}[\cite{RSh}]\label{dimension}
Let be a restricted Lie algebra. Then, for every $m,n\geq 1$, we have
\begin{enumerate}
\item $D_n(L)=\sum_{ip^j\geq n} \gamma_i(L)^{p^j}$,
\item $[D_n(L), D_m(L)]\su \gamma_{m+n}(L)$,
\item $D_n(L)^{p}\subseteq D_{np}(L)$.
\end{enumerate}
\end{thm}

\begin{prop}[\cite{U-IJAC}]\label{dimn-sub-ext}
Let $R$ be a restricted subalgebra of a restricted Lie algebra $L$ and $m$ a positive integer. Then $\w(R)\cap 
\w(L)\w^m(R)=\w^{m+1}(R)$; hence, $L\cap \w(L)\w^m(R)=  D_{m+1}(R)$.
\end{prop}

\begin{thm}[\cite{U-JA08}]
Let $L$ be a restricted Lie algebra and $S$ a restricted Lie subalgebra of $L$. For every positive integer $n$,
the following subalgebras of $L$ coincide.
\begin{enumerate}
\item  $D_{n+2}(S)+D_{n+1}(S\cap D_2(L))$,
\item  $L\cap (\w^{n+2}(S)+\w(S\cap D_2(L))\w^n(S))$,
\item  $L\cap \w^2(L)\w^n(S)$.
\end{enumerate}
\end{thm}

There is a close relationship between restricted Lie algebras and finite $p$-groups.
Indeed, a variant of PBW Theorem was proved by Jennings in \cite{Je} and later 
extended in \cite{U-JPAA}. This analogue of PBW Theorem for group algebras proves to be a very 
useful tool as, for example, one can prove the following Fox-type results.
Below, $\w(G)$ denotes the augmentation ideal of the group algebra $\F G$ over a field $\F$ of positive characteristic $p$. 

\begin{thm}[\cite{U-JPAA}]  Let $G$ be a finite $p$-group. For every subgroup $S$
of $G$ and every positive integer $n$, we have
$$
G\cap (1+\w(G)\w^n(S))=D_{n+1}(S).
$$
\end{thm}

\begin{thm}[\cite{U-JPAA}]  Let $G$ be a finite $p$-group. For every subgroup $S$
of $G$ and every positive integer $n$, we have
$$
G\cap (1+\w^2(G)\w^n(S))=D_{n+2}(S)D_{n+1}(S\cap D_2(G)).
$$
\end{thm}

\section{Powers of the augmentation ideal}\label{aug}
Let $L$ be a Lie algebra with universal enveloping algebra $U(L)$.
A first natural question is whether $U(L)$ determines $\w(L)$. The following lemma answers this
question in the affirmative.

\begin{lem}[\cite{RU}]
Let $L$ and $H$ be Lie algebras and suppose that $\varphi: U(L)\to
U(H)$ is an algebra isomorphism. Then there exists an algebra
isomorphism $\psi:U(L)\to U(H)$ with the property that
$\psi(\omega(L))=\omega(H)$.
\end{lem}

Henceforth, $\varphi:U(L)\to U(H)$ denotes an algebra isomorphism
that preserves the corresponding augmentation ideals.
Since $\varphi$ preserves $\w(L)$, it also preserves the
filtration of $U(L)$ given by the powers of $\w(L)$:
$$
U(L)=\w^0(L) \supseteq\w^1(L) \supseteq \w^2(L)\supseteq \ldots.
$$
Corresponding to this filtration is the graded associative algebra
$$\gr(U(L))=\oplus_{i\geq 0}\omega^i(L)/\omega^{i+1}(L),$$
where the multiplication in $\gr(U(L))$ is induced by
$$(y_i+\w^{i+1}(L))(z_j+\w^{j+1}(L))=y_iz_j+\w^{i+j+1}(L),$$ for all
$y_i\in \w^i(L)$ and $z_j\in \w^j(L)$. Certainly $\gr(U(L))$ is
determined by $U(L)$.

There is an analogous construction for Lie algebras.  That is, one
can consider the graded Lie algebra of $L$ corresponding to its
lower central series given by
$\gr(L)=\oplus_{i\geq1}\gamma_i(L)/\gamma_{i+1}(L)$.
Note that each quotient $\gamma_i(L)/\gamma_{i+1}(L)$ embeds into 
the corresponding quotient $\omega^i(L)/\omega^{i+1}(L)$. Indeed, this way we get a 
Lie algebra homomorphism from $\gr(L)$ into $\gr(U(L))$ which induces an algebra map from 
$U(\gr(L))$ to $\gr(U(L))$. We have:

\begin{thm}[\cite{RU}]\label{graded-iso} For any Lie algebra $L$, the map $\phi:U(\gr(L))\to \gr(U(L))$
is an isomorphism of graded associative algebras.
\end{thm}

Note that under the isomorphism $\phi$ given in Theorem \ref{graded-iso}, we have
$\phi( L/\gamma_2(L))=\w(L)/\w^2(L)$.  Since $\gr(L)$ as a Lie algebra is generated by $ L/\gamma_2(L)$, we deduce that $\phi(\gr(L))$ is the Lie subalgebra of $\gr(U(L))$ genearted by 
$\w(L)/\w^2(L)$. Hence:
 
\begin{cor}[\cite{RU}]
The graded Lie algebra $\gr(L)$ is determined by $U(L)$.
\end{cor}

\begin{cor}[\cite{RU}]\label{lcs} For each pair of integers $(m,n)$ such that $n\ge m\ge 1$,
the quotient $\gamma_n(L)/\gamma_{m+n}(L)$ is determined by
$U(L)$.
\end{cor}

A useful tool that is used to prove many of the results is as follows. Recall that 
 the height of an element $y\in L$, $\nu(y)$,  is 
 the largest integer $n$ such that $y\in
\gamma_n(L)$ if $n$ exists and is infinite if it does not.

\begin{thm}[\cite{Ri, RU}] Let $L$ be an arbitrary Lie algebra and let
$X=\{\bar{x_i}\}_{i\in {\mathcal I}}$ be a homogeneous basis of
$\gr(L)$. Take a coset representative $x_i$ for each $\bar{x_i}$.  Then the
set of all PBW monomials $x_{i_1}^{a_1}x_{i_2}^{a_2}\cdots x_{i_s}^{a_s}$
with the property that $\sum_{k=1}^s a_k\nu(x_{i_k})=n$
forms an $\F$-basis for $\w^n(L)$ modulo $\w^{n+1}(L)$, for every $n\geq 1$.
\end{thm}

A Lie algebra $L$ is called  residually nilpotent if $\cap_{n\geq
1}\gamma_n(L)=0$; analogously, an associative ideal $I$ of
$U(L)$ is residually nilpotent whenever $\cap_{n\geq 1}I^n=0$.

\begin{thm}[\cite{RU}]
Let $L$ be a Lie algebra. Then $L$ is residually nilpotent as a Lie algebra if and only if 
$\w(L)$ is
residually nilpotent as an associative ideal.
\end{thm}

We can now summarize the invariants of $L$ that are determined by $U(L)$.

\begin{thm}[\cite{RU}]\label{major-invariants} The following statements hold for every Lie algebra $L$ over any field.
\begin{enumerate}
\item Whether or not $L$ is residually nilpotent is determined.
\item Whether or not $L$ is nilpotent is determined.
\item\label{nilpotence-class} If $L$ is nilpotent then the nilpotence class of $L$ is determined.
\item If $L$ is nilpotent then the minimal number
of generators of $L$ is determined.
\item If $L$ is a finitely generated free nilpotent Lie algebra then $L$ is determined.
\item The quotient $L'/L''$ is determined.
\item Whether or not $L'$ is residually nilpotent is determined.
\item\label{solvable} If $L$ is finite-dimensional,  then whether or not $L$ is solvable is
determined.
\end{enumerate}
\end{thm}

Part \eqref{solvable} of Theorem \ref{major-invariants} over a field of characteristic zero was proved in \cite{RU}, however,
 according to \cite{V}
enveloping algebra of a finite-dimensional Lie algebra over any field  can be embedded into a (Jacobson)
radical algebra if and only if $L$ is solvable.

\section{Low dimensional nilpotent Lie algebras}\label{low-dim}

Based on results for simple Lie algebras in \cite{M}, it was shown in \cite{CKL}
that $L$ is determined by $U(L)$ in the case when $L$ is any Lie
algebra of dimension at most three over a field of any
characteristic other than two.

In this section we focus on low dimensional nilpotent Lie algebras. 
A classification of such Lie algebra is well known and
can be found, for instance, in \cite{degraaf}. Since there is a unique isomorphism class of nilpotent Lie  algebras
with dimension 1, and there is a unique such class with dimension 2,
the isomorphism problem  is trivial in these
cases. 

Up to isomorphism, there are two nilpotent Lie algebras with dimension $3$ one of which is abelian and the other is non-abelian. By  Part \eqref{nilpotence-class} of  Theorem \ref{major-invariants}, their universal enveloping algebras
must be non-isomorphic. The number of 4-dimensional nilpotent Lie algebras
is 3. One of these algebras is abelian, the second has nilpotency class 2, and
the third has nilpotency class 3. Again, by Part \eqref{nilpotence-class} of  
Theorem \ref{major-invariants},
their universal enveloping algebras are pairwise non-isomorphic.

A strategy for higher dimensions is as follows, which we have used  for dimensions 5 and 6.
For an arbitrary nilpotent Lie algebra  $L$, we know, by Corollary \ref{lcs}, that  the \emph{nilpotency sequence}  $(\dim\gamma_1(L),\dim  \gamma_2(L),\ldots)$, after omitting the tailing zeros, is determined. 
So, nilpotent Lie algebras  of the same finite dimension can be divided into smaller clusters where members of each cluster have the same nilpotency sequence. The investigation of  the isomorphism problem then reduces to the Lie algebras in the same cluster. 
For example, there are 9 isomorphism classes of nilpotent Lie algebras with dimension 5 which can be found in \cite{degraaf}. 
The nilpotency sequence of a nilpotent Lie algebra of dimension 5 is then  one of $(5)$, $(5,1)$, $(5,2)$,  $(5,2,1)$, $(5,3,1)$,
$(5,3,2,1)$, $(5,3,2,1)$. We can now summerize the results for dimensions 5 and 6 as follows:

\begin{thm}[\cite{SU}]\label{5main}
Let $L$ and $H$ be 5-dimensional  nilpotent Lie algebras 
over a field $\F$. If $\uea{L}\cong\uea{H}$, then one of the followings must hold:
\begin{itemize}
\item[(i)] $L\cong H$;
\item[(ii)] $\charac\,\F=2$ and either  $L$ and $H$ are isomorphic to  Lie algebras  $L_{5,3}$ and $L_{5,5}$ or 
$L_{5,6}$ and $L_{5,7}$ in \cite[Section 5]{degraaf}.
\end{itemize}
\end{thm}

\begin{thm}[\cite{SU}]\label{6main}
Let $L$ and $H$ be 6-dimensional  nilpotent Lie algebras 
over a field $\F$ of characteristic not 2. If $\uea{L}\cong\uea{H}$, then one of the followings must hold:
\begin{itemize}
\item[(i)] $L\cong H$.
\item[(ii)] $\charac\,\F=3$ and $L$ and $H$ are isomorphic to one
of the following pairs of Lie algebras in \cite[Section 5]{degraaf}:
$L_{6,6}$ and $L_{6,11}$;
$L_{6,7}$ and $L_{6,12}$; $L_{6,17}$ and $L_{6,18}$; $L_{6,23}$ and $L_{6,25}$.
\end{itemize}
\end{thm}

At the time when Theorem \ref{6main} was proved a complete list of nilpotent Lie algebras of dimension 6 over a field of characteristic 2 was not available. Recently, this list was obtained in \cite{CdS}.
Since, by Theorem \ref{5main}, $U(L_{5,3})\cong U(L_{5,5})$ over a filed of characteristic 2, then it is evident that setting $L=L_{5,3}\oplus \F$ and $H=L_{5,5}\oplus \F$ provides a pair of non-isomorphic nilpotent Lie algebras of dimension 6 such that $U(L)\cong U(H)$ over a field of characteristic 2.

\section{Positive characteristic and restricted Lie algebras}\label{rest}

The results of Section \ref{low-dim}, show in particular that over a field of characteristic 
2 or 3 there exist non-isomorphic nilpotent Lie algebras $L$ and $H$ such that $U(L)\cong U(H)$, thereby providing counterexamples at least in low dimensions. However, the following example provides counterexamples 
over any field of positive characteristic $p$ and dimension $p+3$.

\begin{eg}[\cite{SU}]\label{ex1}
\emph{Let $A=\F x_0+\cdots +\F x_p$ be an abelian Lie algebra over a field $\F$ of characteristic 
$p$. Consider the Lie algebras  $L=A+\F \lambda+\F \pi$ and $H=A+ \F \lambda + \F z$ with 
relations given by  $[\lambda,x_i]=x_{i-1}$,  $[\pi,x_i]=x_{i-p}$, $[\lambda,\pi]=[z,H]=0$, 
and $x_i=0$  for every $i<0$. Then we have:
\begin{enumerate}
\item  $L$ and $H$ are both metabelian and nilpotent of 
class $p+1$. 
\item The centre of $L$ is spanned by $x_{0}$ while the centre of 
$H$ is spanned by $z$ and $x_{0}$; so,
$L$ and $H$ are not isomorphic. 
\item The Lie homomorphism  $\Phi: L\to U(H)$ defined by  $\Phi_{|A+\F \lambda}=\mbox{id}$, 
$\Phi(\pi)=z+\lambda^p$ can be extended to a Hopf algebra isomorphism $U(L)\to U(H)$.
\end{enumerate}
}
\end{eg}

So, the isomorphism problem for enveloping algebras of nilpotent Lie algebras has a negative solution over any field of positive characteristic. Another counterexample can be given in the class of free Lie algebras based on   \cite[Theorem 28.10]{MZ}.
 Recall that the universal enveloping algebra of the free
Lie algebra $L(X)$ on a set $X$ is the free associative algebra
$A(X)$ on $X$. 

\begin{eg}[\cite{RU}]\label{ex2}
\emph{Let $\F$ be a field of odd characteristic $p$ and let $L(X)$ be
the free Lie algebra on $X=\{x,y,z\}$ over $\F$. Set
$h=x+[y,z]+(\mbox{ad } x)^p(z)\in L(X)$ and put $L=L(X)/\langle
h\rangle$, where $\langle h\rangle$ denotes the ideal generated by
$h$ in $L(X)$. Then we have
\begin{enumerate}
\item $L$ is not a free Lie algebra. 
\item There exists a Hopf algebra isomorphism between $U(L)$ and the
2-generator free associative algebra.
\item The minimal number of generators required to generate $L$ is 3.
\end{enumerate}}
\end{eg}

When the underlying field has positive characteristic, it seems natural to  consider the isomorphism problem for   restricted Lie algebras, instead.
\subsection{Restricted isomorphism problem}
 Let $L$ be a restricted Lie algebra with the restricted enveloping algebra $u(L)$ over a field $\F$ of positive characteristic $p$. Let $\w(L)$ denote the augmentation ideal of $u(L)$ which is the kernel of the augmentation map $\epsilon_{_L}: u(L)\to \F$  induced by $x\mapsto 0$, for every $x\in L$.
Let $H$ be another restricted Lie algebra such that
 $\varphi: u(L)\to u(H)$ is an algebra isomorphism.
We observe that the map $\eta: L\to u(H)$ defined by
$\eta=\varphi-\varepsilon _{_H}\varphi$ is a restricted Lie algebra homomorphism.
Hence,  $\eta$ extends to an algebra homomorphism
$\overline{\eta}: u(L)\to u(H)$. In fact,
$\overline{\eta}$ is an isomorphism that preserves the augmentation ideals, that is
$\overline{\eta}(\w(L))=\w(H)$. So, without loss of generality, we assume that $\varphi:u(L)\to u(H)$ is an algebra isomorphism that preserves the augmentation ideals. 

Note that the role of lower central series in Lie algebras is played by the dimension subalgebras in restricted Lie algebras. Recall from Theorem \ref{dimension} that  the $n$-th dimension subalgebra of $L$ is 
$$
D_n(L)=L\cap \w^n(L)=\sum_{ip^j\geq n} \gamma_i(L)^{p^j}.
$$

Now,  consider the graded restricted Lie algebra:
$$
\gr(L):=\bigoplus_{i\ge 1} D_i(L)/D_{i+1}(L),
$$
where the Lie bracket and the $p$-map are defined over homogeneous elements and then extended linearley: 
\begin{align*}
[x_i+D_{i+1}(L), x_j+D_{j+1}(L)]&=[x_i, x_j]+D_{i+j+1}(L),\\
(x_i+D_{i+1}(L))^{[p]}&=x_i^p+D_{ip+1}(L)
\end{align*}
 for all
$x_i\in D_i(L)$ and $x_j\in D_j(L)$.
In close analogy with Theorem \ref{graded-iso},  one can see that  $u(\gr(L))\cong \gr(u(L))$ as algebras.
So we may identify $\gr(L)$ as the graded restricted Lie subalgebra of $\gr(u(L))$
generated by $\w^1(L)/\w^2(L)$. Thus, $\gr(L)$ is determined.

Recall that $L$ is said to be in the class $\cFp$ if $L$ is finite-dimensional and $p$-nilpotent.
 Whether or not $L\in \cFp$ is determined by the following lemma, see \cite{RSh}.

\begin{lem}\label{w(L)-nilpotent}
Let $L$ be a restricted Lie algebra. Then $L\in \cFp$ if and only if $\w(L)$ is nilpotent.
\end{lem}
\begin{lem}[\cite{U-PAMS}]\label{nilp-class}
 If $u(L)\cong u(H)$ then the following statements hold.
 \begin{enumerate}
 \item  If $L\in \cFp$ then $\mid cl(L)-cl(H)\mid \leq 1$.
 \item $D_i(L)/D_{i+1}(L)\cong D_i(H)/D_{i+1}(H)$, for every $i\geq 1$.
 \end{enumerate}
 \end{lem}

  We remark that
methods of  \cite{RSh} and \cite{RU}  can be adapted to prove that
the quotients $D_n(L)/D_{2n+1}(L)$ and $D_{n}(L)/D_{n+2}(L)$ are also determined, for every $n\geq 1$.
In particular, $L/D_3(L)$ is determined.

Unlike the isomorphism problem in which abelian Lie algebras are determined by their enveloping algebras, the abelian case for the restricted isomorphism problem is not trivial.
Note that if  $L$ is an abelian restricted Lie algebra then the $p$-map reduces to
$$
(x+y)^p=x^p+y^p, \quad (\alpha x)^p=\alpha^px^p,
$$
for every $x,y\in L$ and $\alpha\in \F$. Thus the $p$-map is a semi-linear transformation.

\begin{thm}[\cite{U-PAMS}]\label{prop-perfect}
Let $L\in \cFp$ be an abelian  restricted Lie algebra over a perfect field  $\F$. If $H$ is  a restricted Lie
algebra such
that $u(L)\cong u(H)$, then $L\cong H$.
\end{thm}

\begin{cor}\label{L/L'_p}
Let  $L\in \cFp$ be a restricted Lie algebra over a perfect field. Then $L/L'_p$ is  determined.
\end{cor}

It turns out that over an  algebraically closed field stronger results hold. 

\begin{thm}[\cite{U-PAMS}]\label{prop-alg-closed}
Let $L$ be a finite-dimensional abelian  restricted Lie algebra over an algebraically closed field $\F$. Let $H$
be a restricted Lie algebra such that $u(L)\cong u(H)$. Then $L\cong H$.
\end{thm}

Using the identity $[ab,c]=a[b,c]+[a,c]b$ which holds in any associative algebra, we can see that
 $L'_pu(L)=[\w(L),\w(L)]u(L)$. Thus the ideal $L'_pu(L)$ is preserved by $\varphi$.
Now write $J_L=\w(L)L'+L'\w(L)=\w(L)L'_p+L'_p\w(L)$. Since both $\w(L)L'_p$ and $L'_p\w(L)$ are determined, it
follows that $J_L$ is determined.

\begin{thm}[\cite{U-PJM}]\label{dim-L-mod}
If $L\in \cFp$ and $\F$ is perfect  then $L/(L'^p{+}\gamma_3(L))$ is determined.
\end{thm}

\begin{thm}[\cite{U-PAMS}]\label{L'-quo}
Suppose that $L$ and $H$ are finite-dimensional restricted Lie algebras such that $u(L)\cong u(H)$. Then,
for every positive integer $n$, we have
$$
D_n(L'_p)/D_{n+1}(L'_p)\cong D_n(H'_p)/D_{n+1}(H'_p).
$$ 
\end{thm}

\begin{lem}[\cite{U-PAMS}]\label{exp-H'}
Let $L\in \cFp$  such that $cl(L)=2$. Then, $\dim_{\F} {L'_p}^{p^t}$ is determined, for every $t\geq 0$.
\end{lem}

\begin{lem}[\cite{U-PAMS}]\label{dim-H'}
Let $L\in \cFp$   such that $L'_p$ is cyclic. The following statements hold.
\begin{enumerate}
\item $cl(L)\leq 3$.
\item  We have $L'^{p^t}u(L)=(L'_pu(L))^{p^t}$,  for every $t\geq 1$.
\end{enumerate}
\end{lem}

A restricted Lie algebra $L$ is called \emph{metacyclic} if $L$ has a cyclic restricted  ideal $I$ such that 
$L/I$ is cyclic.  Recall that a $p$-polynomial in $x$ has the form $c_0x+c_1x^p+\cdots+c_t x^{p^t}$, where each $c_i\in \F$. So, if $L$ is metacyclic then  there exist generators $x,y\in L$ and  some
$p$-polynomials $g$ and $h$ such that 
$$
h(x)\in \la y\ra_p,\quad [y,x]=g(y).
$$
Now let $L$ be a non-abelian  metacyclic restricted Lie algebra in the class $\cFp$.
It turns out  that there exist another $p$-polynomial $f$ and positive integers $m,n$ such that the following
 relations hold in $L$:
\begin{align*}
& x^{p^m}=f(y)=y^{p^r}+\cdots,\\
& y^{p^n}=0,\\
& [y,x]=g(y)=b_sy^{p^s}+\cdots, b_s\neq 0.
\end{align*}
Since $L$ is not abelian, we have $1\leq r\leq n$ and $1\leq s\leq n-1$.

\begin{thm}[\cite{U-PAMS}]
 Let $L\in \cFp$ be a metacyclic restricted Lie algebra over a perfect field of positive characteristic. Then $L$
 is determined by $u(L)$.
 \end{thm}

\section{Other observations}\label{hopf}

Because enveloping algebras are Hopf algebras, it also makes sense
to consider an enriched form of the isomorphism problem that takes
this Hopf structure into account.

Recall that a bialgebra is a vector space ${\mathcal H}$ over a
field $\F$ endowed with an algebra structure $({\mathcal H},
M, u)$ and a coalgebra structure
$({\mathcal H}, \Delta, \epsilon)$ such
that $\Delta$ and $\epsilon$ are algebra
homomorphisms. A bialgebra ${\mathcal H}$ having an antipode
$S$ is called a Hopf algebra. It is well-known that
the  enveloping algebra of a (resticted) Lie algebra is a  Hopf
algebra, see for example \cite{BMPZ} or \cite{MZ}. Indeed, the
counit $\epsilon$ is the augmentation map and the coproduct $\Delta$ is induced by
$x\mapsto x\otimes 1+1\otimes x$, for every $x\in L$. An explicit
description of  $\Delta$ can be given in terms of a PBW
monomials (see, for example, Lemma 5.1 in Section 2 of
\cite{SF}). The antipode $S$ is induced by $x\mapsto -x$,
for every $x\in L$. The following proposition is well-known
(see Theorems 2.10 and 2.11 in Chapter 3 of \cite{BMPZ}, for example).

\begin{prop}
Let $L$  be a Lie algebra over a field $\F$ of characteristic $p\ge0$.
\begin{enumerate}
\item If $p=0$ then the set of primitive elements of $U(L)$ is $L$. Since an isomorphism of 
Hopf algebras perserves the primitive elements, the Hopf algebra structure of $U(L)$ determines $L$.
\item If $p>0$ and $L$ is a restricted Lie algebra then the primitive elements of $u(L)$ is $L$. Hence, 
the Hopf algebra structure of $u(L)$ determines $L$.
\end{enumerate}
\end{prop}

 If $p>0$ then the set of primitive elements of $U(L)$ is $L_p$,
the restricted Lie subalgebra of $U(L)$ generated by $L$. Thus,
any Hopf algebra isomorphism from $U(L)\to U(H)$
restricts to a restricted Lie algebra isomorphism $L_p\to H_p$.

We now present an example illustrating that the
analogous isomorphism problem for enveloping algebras of Lie
superalgebras fails utterly in the sense that  $\dim_{\F} L$ may not be determined.

Let $\F$ be a field of characteristic not 2. In the case of characteristic 3, we add the
axiom $[x,x,x]=0$ in order for the universal enveloping algebra,
$U(L)$, of a Lie superalgebra $L$ to be well-defined.

\begin{eg}[\cite{RU}]
\emph{Let $L=\F x_0$ be the free Lie superalgebra on one generator $x_0$
of even degree, and let $H=\F x_1+\F y_0$ be the free Lie
superalgebra on one generator $x_1$ of odd degree, where $y_0=[x_1,x_1]$.
Then $U(L)$ is
isomorphic to the polynomial algebra $\F[x_0]$ in the indeterminate
$x_0$. On the other hand, $U(H)\cong \F[x_1,y_0]/I$, where
$I$ is the ideal of the polynomial algebra $\F[x_1,y_0]$ generated by
$y_0-2x_1^2$. Hence, $U(H)\cong \F[x_1]\cong \F[x_0]\cong U(L)$.
However, $L$ and $H$ are not isomorphic  since they do not  have the same
dimension.}
\end{eg}

\section{Open Problems}\label{quo}
Below we list a set of problems that are interesting to investigate:
\begin{enumerate}
\item An interesting open problem asks whether or not similar examples as Example \ref{ex1} can
occur in characteristic zero; that is, does there exist a non-free
Lie algebra $L$ over a field of characteristic zero such that
$U(L)$ is a free associative algebra? 
\item Is the derived length of a solvable Lie algebra determined?
\item Let $L$ be a finite-dimensional Lie algebra over a field of characteristic zero. Is $Z(L)$ determined?
\item Let $L$  be a finite-dimensional metabelian Lie algebra over a field of characteristic zero. Is $L$ determined?
\item \textbf{Conjecture:}  Let $L$  be a finite-dimensional nilpotent Lie algebra over a field of characteristic zero. Then $L$ is determined by $U(L)$.
\item Provide a counterexample to the restricted  isomorphism problem.
\item Let $L$ be a finite-dimensional restricted  Lie algebra over a field of positive characteristic. Is $Z(L)$ determined?
\item Let $L$ be a finite-dimensional $p$-nilpotent restricted Lie algebra over a field of positive characteristic. Is the nilpotence class of $L$ determined?
\item Let $L$ be a finite-dimensional $p$-nilpotent restricted Lie algebra over a perfect field of positive characteristic $p$ . Is  $L$ determined by $u(L)$?
\end{enumerate}


\begin{thebibliography}{99}\label{bib}



\bibitem{B} Y.A. Bahturin, {\em Identical Relations in Lie Algebras}
(VNU Science Press, b.v., Utrecht, 1987).

\bibitem{BMPZ} Y.A. Bahturin, A. Mikhalev, V. Petrogradsky, M. Zaicev, {\em 
Infinite-Dimensional Lie Superalgebras}, de Gruyter Exp. Math. {\bf 7} (de Gruyter, Berlin, 
1992).

\bibitem{CdS}
S. Cical\`o, W.A.  de Graaf, C.  Schneider, 
Six-dimensional nilpotent Lie algebras, 
\emph{Linear Algebra Appl.} \textbf{436} (2012), no. 1, 163--189. 

\bibitem{CKL} J. Chun, T. Kajiwara, J. Lee, Isomorphism theorem on low dimensional Lie 
algebras, {\em Pacific J. Math.} {\bf 214} (2004), no. 1, 17--21.





\bibitem{degraaf} W.A. de Graaf,
\newblock
Classification of 6-dimensional nilpotent Lie algebras over fields of characteristic not 2,
\newblock {\em J. Algebra} {\bf 309} (2007), no. 2, 640--653.

\bibitem{Fox} R.H. Fox, Free differential calculus I,
\emph{Ann. of Math.} \textbf{57} (1953), 547--560.


\bibitem{Gup} N.D. Gupta, A problem of R.H. Fox,
{\em Canad. Math. Bull.} {\bf 24} (1981), 129--136.






\bibitem{H} M. Hertweck, A counterexample to the isomorphism problem for integral group 
rings, {\em Ann. of Math. (2)} {\bf 154} (2001), no. 1, 115--138.

\bibitem{Hur2} T.C. Hurley, Identifications in a free group,
\emph{J. Pure and Applied Algebra} \textbf{48} (1987), 249--261.

\bibitem{HSe} T.C. Hurley, S.K. Sehgal, Groups related to Fox subgroups,
\emph{Comm. Algebra} \textbf{28} (2000), no. 2, 1051--1059.
 
\bibitem{Je} S.A. Jennings, The group ring of a class of infinite nilpotent groups,
\emph{Canad. J. Math.} \textbf{7} (1955), 169--187.



\bibitem{M} P. Malcolmson, Enveloping algebras of simple three-dimensional Lie algebras,
{\em J. Algebra} {\bf 146} (1992), 210-218.

\bibitem{MZ} A.A. Mikhalev, A.A. Zolotykh, \emph{Combinatorial aspects of Lie superalgebras}
(CRC Press, Boca Raton, FL, 1995).


\bibitem{Ri} D.M. Riley, The dimension subalgebra problem for enveloping algebras of Lie
superalgebras, {\em Proc. Amer. Math. Soc.} {\bf 123} (1995), no. 10, 2975--2980.

\bibitem{RSh}  D.M. Riley, A. Shalev, Restricted Lie algebras and their envelopes,
 \emph{Canad. J. Math.} \textbf{47} (1995), 146--164.


\bibitem{RU} D.M. Riley, H. Usefi, The isomorphism problem for enveloping algebras,
{\em Alg. Repr. Theory}, \textbf{10} (2007), no. 6, 517--532.

\bibitem{RS} K. Roggenkamp, L. Scott, Isomorphisms of $p$-adic group rings,
{\em  Ann. of Math. (2)} {\bf 126} (1987), no. 3, 593--647.

\bibitem{SF} H. Strade, R. Farnsteiner, {\em  Modular Lie Algebras and Their Representations}, 
Monographs and Textbooks in Pure and Applied Mathematics {\bf 116} (Dekker, New York, 1988).

 \bibitem{SU} C. Schneider,  H. Usefi, Isomorphism problem for enveloping algebras of nilpotent Lie algebras,  \emph{Journal of Algebra}, {\bf 337} (2011), 126-140.


 \bibitem{U-PJM} H. Usefi, Isomorphism invariants of restricted enveloping algebras, \emph{Pacific Journal of Mathematics},  \textbf{246} (2010), No. 2, 487-494.
 
 \bibitem{U-PAMS} H. Usefi, The restricted isomorphism problem for  metacyclic restricted Lie algebras, \textit{Proceedings of the American Mathematical Society}, \textbf{136} (2008), 4125-4133.

 \bibitem{U-JA08} H. Usefi, Fox-type problems in enveloping algebras, \textit{Journal of Algebra}, \textbf{319} (2008) 2489-2495.


 \bibitem{U-JPAA} H. Usefi,  Identifications in modular group algebras,
 \emph{Journal of Pure and Applied Algebra}, \textbf{212} (2008) 2182--2189.

 \bibitem{U-IJAC} H. Usefi, The Fox problem for   free restricted Lie  algebras, \emph{International Journal of Algebra and Computation}, \textbf{18} (2008), no. 2, 271-283. 

\bibitem{V} A.I. Valitskas, A representation of finite-dimensional Lie algebras in radical 
rings, {\em Dol. Akad. Nauk SSSR} {\bf 279} (1984), no. 6, 1297-1300.

\bibitem{W} A. Weiss, Rigidity of $p$-adic $p$-torsion,
{\em Ann. of Math. (2)} {\bf 127} (1988), no. 2, 317--332.















\bibitem{Yun} I.A. Yunus, A problem of Fox,
{\em Dokl. Akad. Nauk SSSR} {\bf 278} (1984), no. 1, 53--56.

\bibitem{Yun2} I.A. Yunus,  The Fox problem for Lie algebras,
 \emph{Uspekhi Mat. Nauk} \textbf{39} (1984), no. 3,  251--252.



\end{thebibliography}
\end{document}